\documentclass[11pt]{article} 
\usepackage{amsmath, amsfonts, amssymb}
\usepackage{hyperref}
\newcounter{Mycounter}[section]
\newcounter{lemma}[section]
\newcounter{claim}[section]
\newcounter{sublemma}[section]
\newcounter{corollary}[section]
\newcounter{theorem}[section]
\renewcommand{\thetheorem}{{Theorem \thesection.\arabic{theorem}}}
\newcommand{\theorem}{%
     \setcounter{theorem}{\value{Mycounter}}
     \refstepcounter{theorem}
     \stepcounter{Mycounter}
     {\bf \thetheorem:\ }}

\newcounter{conjecture}[section]
\renewcommand{\theconjecture}{{Conjecture \thesection.\arabic{conjecture}}}
\newcommand{\conjecture}{%
     \setcounter{conjecture}{\value{Mycounter}}
     \refstepcounter{conjecture}
     \stepcounter{Mycounter}
     {\bf \theconjecture:\ }}

\newcounter{proposition}[section]
\newcounter{definition}[section]
\renewcommand{\thedefinition}
       {{Definition~\thesection.\arabic{definition}}}
\newcommand{\definition}{%
     \setcounter{definition}{\value{Mycounter}}
     \refstepcounter{definition}
     \stepcounter{Mycounter}
     {\bf \thedefinition:\ }}
     
\newenvironment{proof}{
{\bf Proof:\\}}
{
\begin{flushright}
•$\blacksquare$
\end{flushright}}

\begin{document}
\begin{center}

{\LARGE\bf The second Betti number of hyperk\"ahler manifolds}
\\[2mm]
Nikon Kurnosov\footnote{ The author is partially
supported by the Simons Foundation grant and AG Laboratory 
NRU-HSE, RF government grant, ag. 11.G34.31.0023. 
}
\\[2mm]
{\it nikon.kurnosov@gmail.com}
\end{center}

{\small 
\hspace{0.15\linewidth}
\begin{minipage}[t]{0.75\linewidth}
{\bf Abstract} \\
Let $M$ be a compact irreducible hyperk\"ahler manifold, from Bogomolov inequality \cite{V1} we obtain forbidden values of the second Betti number $b_2$ in arbitrary dimension.\\
\textbf{UPD:} Unfortunately, decomposition of dual to BBF-form is not right in the main theorem. Instead of this work, take a look on recent preprints of Sawon and me on boundedness of $b_2$ for hyperk\"ahler manifolds.
\end{minipage}
}
{
\small
\tableofcontents
}

\section{Introduction}

One of the main conjectures in the theory of hyperk\"ahler manifolds is finitness of number of classes of hyperk\"ahler manifolds up to deformation in each dimension \cite{Bea1}. Huybrechts \cite{H1} has proved finiteness of the number of deformation classes
of holomorphic symplectic structures on each smooth manifold. Moreover, in most
dimensions as we see below we only know the two examples due to Beauville
\cite{Bea1}. Given the properties of a hyperk\"ahler
structure, one might be led to suspect that there can not be many more. In two
dimensions there are two more examples constructed by O'Grady \cite{O1, O2}. The conjecture
is that there are finitely many hyperk\"ahler manifolds up to
the deformation in each dimension. It is well-known, e.g. \cite{H2}, for compact
hyperk\"ahler manifolds the second Betti number $\geqslant$ 3. Thanks to Guan \cite{G} it is known for 4-dimensional hyperk\"ahler manifolds that $ 3\leq b_2 \leq 8$ or $ b_2 =23$ and also there are some bounds on $ b_3$. The conjecture is that the second Betti number is bounded in any dimension. In the present paper, we obtain results about the Betti numbers of compact hyperk\"ahler manifolds and state the following

\hfill

{\bf Theorem:} Let $M$ be a compact irreducible hyperk\"ahler manifold of complex dimension $n$. Then 
\begin{center}
$b_2 \neq \frac{10 + n^2 - n}{2}$,
\end{center}
where $b_2$ is the second Betti number of $M$.
\hfill

\section{Hyperk\"ahler manifolds}
\hfill
\definition \label{Def_HK_manifolds} 
(\cite{Bes}) Let $(M, g)$ be a Riemannian manifold, and $I$, $J$,$K$ are
endomorphisms of a real tangent bundle satisfying the relation 
$I\circ J=-J\circ I = K$. If the metric $g$ on $M$ is K\"ahler with respect to these complex 
structures than $M$ is called a {\bf hyperk\"ahler manifold}.

\hfill

Clearly, complex dimension of hyperk\"ahler manifold $M$ is even. A compact K\"ahler manifold is hyperk\"ahler (HK) if it is simply connected and the space of its global holomorphic two-forms is spanned by a symplectic form. In algebraic geometry the word "hyperk\"ahler" is synonymous with "holomorphic symplectic". 

\hfill

\definition \label{Def_holom_symplectic} 
A manifold $M$ is called {\bf holomorphically symplectic} if it is a complex manifold with a closed holomorphic 2-form $\Omega$ over $M$ such that 
$\Omega^n = \Omega \wedge \Omega \wedge ... \wedge \Omega$ is
a nowhere degenerate section of a canonical class of $M$, where $2n = \dim_\mathbb{C}(M)$.

\hfill

Consider the K\"ahler forms $w_I(\cdot, \cdot) := g(\cdot, I\cdot), w_J(\cdot, \cdot) := g(\cdot, J\cdot), w_K(\cdot, \cdot) := g(\cdot, K\cdot)$ on $M$.
A simple algebraic calculation \cite{Bes} shows that the following form
\[
\Omega = w_J+\sqrt{-1}w_K
\]
is of type $(2,0)$ on $(M, I)$. It is closed, holomorphic and moreover nowhere degenerate, as another linear algebraic argument shows, hence it is a holomorphic symplectic form. Thus, the underlying complex manifold
$(M, I)$ is holomorphically symplectic. The converse assertion
is also true:

\hfill

\theorem \label{symplectic=>HK}
(\cite{Bea2}, \cite{Bes})
Let $M$ be a compact, K\"ahler, holomorphically
symplectic manifold with the 
holomorphic symplectic form $\Omega$, $w$ its K\"ahler form, $n = \dim_\mathbb{C} M$. Assume that
$\int_M w^n = \int_M (Re \Omega)^n$.
Then there is a unique hyperk\"ahler 
structure $(I, J, K, g)$
over $M$ such that the cohomology class of the symplectic form
$w_I = g(\cdot, I\cdot)$ is equal to $w$ and $\Omega = w_J + \sqrt{-1}w_K$.

\hfill

The Bogomolov-Beauville-Fujiki form was defined in \cite{Bo} and \cite{Bea2}, but it
is easiest to describe it using the Fujiki theorem.

\hfill

\theorem \label{Fujiki}
(\cite{F})
Let $M$ be a simple hyperk\"ahler manifold, $\eta \in H^2(M)$, and $2n = \dim M$. Then $\int_M \eta^{2n} = \lambda q(\eta, \eta)^n$, for some primitive integer quadratic
form $q$ on $H^2(M, \mathbb{Z})$, and $\lambda > 0$ an integer number.

\hfill

Fujiki formula (\ref{Fujiki}) determines the form $q$ uniquely up
to a sign. For odd $n$, the sign is unambiguously determined as well. For even $n$, one singles out one of the two choices by imposing the inequality
\[
q(\Omega, \bar{\Omega}) > 0, \qquad 0 \neq \Omega \in H^{2,0}(M)\]

The two-dimension examples of irreducible holomorphic symplectic manifolds
are called K3 surfaces. In higher dimensions there are
only few examples known. Here is the list of known examples, where manifolds of
the same deformation type are not distinguished.

(i) If $X$ is a K3 surface then the Hilbert scheme Hilb$^n(X)$ is an irreducible
holomorphic symplectic manifold \cite{Bea2}. Its dimension is $2n$ and for $n > 1$ its
second Betti number is equal to 23. Construction of Hilbert Scheme can be
found elsewhere, e.g. \cite{Bea1, Bea2}. Let $X$ be a K3 surface. Take
the symmetric product $X^{\left( r \right)} = X^r /\mathfrak{S}_r$ which
parametrizes subsets of $r$ points in K3 surface $X$, counted with multiplicities;
it is smooth on the open subset $X_0$ consisting of subsets with $r$ distinct points,
but singular otherwise. We obtain a smooth compact manifold, which is called
the Hilbert scheme $X^{\left[ r \right]}$ if we replace ``subset'' by
``subspace''. The natural map $X^{\left[ r \right]} \rightarrow X^{\left( r
\right)}$ is an isomorphism above $X_0$, but it resolves the singularities of $X^{\left( r \right)}$.

Let us describe the easiest case Hilb$^2(X)$ explicitly. For any surface $X$ the
Hilbert scheme Hilb$^2(X)$ is the blow-up Hilb$^2 (X) \rightarrow S^2
(X)$ of the diagonal $\Delta =\{\{x, x\} \left| \right. x \in X\} \subset S^2
(X) =\{\{x, y\} \left| \right. x, y \in X\}$. Equivalently, Hilb$^2(X)$ is the
$\mathbb{Z}/ 2\mathbb{Z}$-quotient of the blow-up of the diagonal in $X
\times X$. Since for a K3 surface there exists only one $\mathbb{Z}/
2\mathbb{Z}$-invariant two-form on $X \times X$, the holomorphic symplectic
structure on Hilb$^2(X)$ is unique.

(ii) If $X$ is a complex torus of dimension two, then the generalized Kummer
variety K$_n(X)$ is an irreducible holomorphic symplectic manifold \cite{Bea2}. Its
dimension is $2n$ and for $n > 2$ its second Betti number is 7. \ The Hilbert
scheme $X^{\left[ n \right]}$ of two-dimensional torus has the same properties
as $X^{\left[ r \right]}$, but it is not simply connected. This is fixed by
considering the composition of maps: $X^{\left[ n \right]} \rightarrow
X^{\left( n \right)} \rightarrow X$, where the last map is $s \left( t_1,
\ldots, t_n \right) = t_1 + \ldots + t_n$. The fibre K$_n(X)$ is a hyperk\"ahler
manifold of dimension $2n$ \ - generalized Kummer manifold.

(iii) O'Grady's 10-dimensional example \cite{O1}. Let again $S$ be a K3 surface, and
M the moduli space of stable rank 2 vector bundles on $S$, with Chern classes
$c_1 = 0, c_2 = 4$. It admits a natural compactification $M$, obtained by adding classes of semi-stable torsion free sheaves. It is singular
along the boundary, but O'Grady \cite{O1} constructs a desingularization of $M$ which
is a new hyperk\"ahler manifold, of dimension 10. Its second Betti number is 24
\cite{R}. Originally, it was proved that it is at least 24 \cite{O1}.

(iv) O'Grady's 6-dimensional example \cite{O2}. The similar construction can be
done starting from rank 2 bundles with $c_1 = 0, c_2 = 2$ on a 2-dimensional
complex torus, that gives a new hyperk\"ahler manifold of dimension 6 as in (iii). Its second Betti number is 8.

Thus we have two series, (i) and (ii), and two sporadic examples, (iii) and
(iv). All of them have different second Betti numbers. It has been proved \cite{KLS} that the moduli spaces for all other sets of numerical parameters unless Hilb$^n(K3)$ do not admit a smooth symplectic resolution of singularities. Note that in any given dimension
and for any given second Betti number $b_2$ one knows at most one real manifold
carrying the structure of an irreducible holomorphic symplectic manifold. So,

\hfill

\conjecture (\cite{Bea1, S}) \label{finitness}
The number of deformation types
of compact irreducible hyperk\"ahler is finite in any dimension (at least for given $b_2$).
\hfill

\section{On the Betti numbers of hyperk\"ahler manifolds}
Let $M$ be a compact K\"ahler manifold of complex dimension $n$. The Hodge number
$h^{p,q}$ denotes the dimension of the corresponding Dolbeault cohomology
\[ h^{p, q} = h^{n - p, n - q} = h^{q, p} \]
If $M$ be a compact connected hyperk\"ahler manifold of real dimension $4m$. We can
find out more equalities on Hodge numbers. By studying the action of $Sp(m)$ on
spaces of harmonic forms, Wakakuwa \cite{W} proved that $b_{2 k} \geqslant
\left(\begin{array}{c}
  k + 2\\
  2
\end{array}\right)$ for $k \leqslant m$ and that the odd Betti numbers $b_{2 k
+ 1}$ of $M$ are all divisible by 4. Moreover, Fujiki using Hodge decompositions
refined relative to a choice of complex structure. He proved \cite{F} that $h^{p,
q} \geqslant h^{p + 1}, q - 1$ whenever $p \geqslant q$.

{\bf Note:} The first Betti number $b_1$ for compact irreducible hyperk\"ahler manifolds is always zero, but
in general odd Betti number not, e.g. $b_3 \left( K_2 \left( T \right) \right)
= 8$.

Wedging with the holomorphic symplectic form $\Omega = w_J + \sqrt{- 1}
w_K$ induces a mapping $H^{p, q} \rightarrow H^{p + 2, q}$ which is
injective for $p + 1 \leqslant m$ and its $(m-p)$-fold iteration is an
isomorphism. In this way, the equalities on Hodge numbers above are
supplemented by the equations
\[ h^{p, q} = h^{2 m - p, q}, 0 \leqslant p, q \leqslant 2 m. \]
Proof of these results has been given by Verbitsky \cite{V2} by considering the
action of the Lie algebra $\mathfrak{s}\mathfrak{o} \left( 5 \right)$ on
cohomology.

Salamon proved \cite{Sa} that if $X$ is a compact hyperk\"ahler manifold of dimension
$2n (=4m)$ then
\[ \sum_{i = 0}^{4 m} \left( - 1 \right)^i \left( 6 i^2 - n \left( 6 n + 1
   \right) b_i \right) = 0 \]
Using Wakakuwa results this shows in particular that $n e (X)$ is divisible
by 24, where $e(X)$ is the Euler characteristic. Note that for a K3 surface the Euler number is 24. Indeed, we can
rewrite Salamon's relation above in terms of Euler characteristics:
\[ \sum_{i = 0}^{4 m} \left( - 1 \right)^i 6 i^2 b_i = n \left( 6 n + 1
   \right) e \left( X \right) . \]

The case of 4-dimensional hyperk\"ahler manifolds has been discussed by Guan \cite{G}.

\hfill

\begin{theorem} \label{Guan's results}
  If $M$ is an irreducible compact hyperk\"ahler manifold of complex dimension 4,
  then
  \begin{itemize}
  
     \item if $b_2 = 23$, then $b_3 = 0$. The Hodge diamond of $M$ is the same as
    that of the Hilbert scheme of pairs of points on a K3 surface.
    
    \item if $b_2 \neq 23$, then $b_2 \leqslant 8$, and if $b_2 = 8$, then
    $b_3 = 0$.
    
    \item in the case of $b_2$ = 7, $b_3$ = 0 or 8.
    
    \item in the case of $b_2 = 3, 4, 5, 6$, then the following cases
    are possible
    \begin{center}
   \begin{tabular}{|l|l|l|l|l|}
  \hline
  $b_2$ & 3 & 4 & 5 & 6\\
  \hline
  $b_3$ & $4 l, l \leqslant 17$ & $4 l, l \leqslant 15$ & $4 l, l \leqslant 9$
  & $4 l, l \leqslant 4$\\
  \hline
\end{tabular}.
    \end{center}
        
    \item the second Chern class lies in the algebra $H^{\left( 4 \right)}$
    generated by $H^2 \left( M \right)$ iff 
    \[ \left( b_2, b_3 \right) = (5,
    36), (7, 8), (8, 0), (23, 0).\]
  \end{itemize}
\end{theorem}

To prove the \ref{Main Theorem} stated in the Introduction we will use Bogomolov inequality \cite{V1}. Denote by $w$ the K\"ahler form on $M$ and let $n$ be $\dim_\mathbb{C} M$, then by Gauss-Bonnet formula,
the cohomology class of $Tr(\Theta\wedge \Theta)$ can be expressed
via $c_1(F)$, $c_2(F)$:
\[ 
  \frac{\sqrt{-1}}{{2\pi}^2} Tr(\Theta\wedge \Theta) = 2 c_2(F) - \frac{n-1}{n} c^2_1(F),
\]
where $F$ is torsion-free coherent sheaf and $\theta$ is curvature of hyperholomorphic connection on it.
Therefore, the following integral (since $c_1=0$)
\begin{equation}
\int_Mc_2\wedge w^{n-2} > 0
\end{equation}
is positive that gives Bogomolov inequality \cite{V1}.

\hfill

\theorem \label{Main Theorem}  Let $M$ be a compact irreducible hyperk\"ahler manifold of complex dimension $n$. \begin{center}
$b_2 \neq \frac{10 + n^2 - n}{2}$,
\end{center}
where $b_2$ is the second Betti number of $M$.
\\
\hfill

\begin{proof}
Let ${e_1,...,e_{b_2}}$ be an orthonormal basis in $H^2(M)$ with respect to the Bogomolov-Beauville quadratic form $q$. Denote the dual form of $q$ by $Q \in $ Sym$^2H^2(M)$. Then
\[ Q = e_1^2 + e_2^2 + e_3^2 - \sum_4^{b_2} e_i^2. \]
Since $c_2 = \mu Q + p$, where $\mu \in \mathbb{Q}$ and $p \in $ (Sym$^2H^2(M))^\perp$ then
\[ \int_M Q \wedge e_1^{n-2} \neq 0, \]
otherwise Bogomolov inequality (1) is not fulfilled.

By Fujiki formula we obtain
\[ \int_M e_1^2 \wedge e_1^{n-2} = n!,\]
\[ \int_M e_i^2 \wedge e_1^{n-2} = 2 \cdot (n - 2)! \qquad (i \neq 1).\]
Then using composition of $Q$ given above and sum up
\[ Q \wedge e_1^{n - 2} = n! + 2 \cdot \left( n - 2 \right) ! + 2 \cdot
\left( n - 2 \right) ! - 2 \cdot \left( n - 2 \right) ! \cdot \left( b_2 - 3
\right). \]
After calculations we obtain
\[ Q \wedge e_1^{n - 2} = \left( n - 2 \right) ! \cdot \left( 10 + n \left( n - 1 \right)
\right) - 2 \cdot \left( n - 2 \right) ! \cdot b_2. \]
Thus,
\[10 + n^2 - n \neq 2 b_2. \]
\end{proof}

In the case of 4-dimensional hyperk\"ahler manifolds $ b_2\neq11$ that has been also proved by Guan (\ref{Guan's results}) earlier and in the dimension 6 impossible value for $ b_2$ is 20. Recall that all known examples of hyperk\"ahler manifolds in dimension 6 have $b_2$ equal to 7 (O'Grady example), 8 (generalized Kummer manifold K$_3(X)$) or 23 (Hilb$_3(X)$).
\\[1.0pt]
\subsection*{Acknowledgements} Author would like to thank M. Verbitsky for valuable discussions and ideas.

\hfill

\noindent {\sc Nikon Kurnosov\\
Laboratory of Algebraic Geometry and its applications,\\
National Research University Higher School of Economics\\
7 Vavilova Str., Moscow, Russia, 117312;\\
Independent University of Moscow,\\
11 Bol.Vlas'evskiy per., Moscow, Russia, 119002}\\
\it  nikon.kurnosov@gmail.com

\end{document}